# EXTENSION OF INCOMPRESSIBLE SURFACES ON THE BOUNDARY OF 3-MANIFOLDS

Michael Freedman, Hugh Howards and Ying-Qing Wu

ABSTRACT. An incompressible surface $F$ on the boundary of a compact orientable 3-manifold $M$ is arc-extendible if there is an arc $\gamma$ on $\partial M - \text{Int} F$ such that $F \cup N(\gamma)$ is incompressible, where $N(\gamma)$ is a regular neighborhood of $\gamma$ in $\partial M$. Suppose for simplicity that $M$ is irreducible, and $F$ has no disk components. If $M$ is a product $F \times I$, or if $\partial M - F$ is a set of annuli, then clearly $F$ is not arc-extendible. The main theorem of this paper shows that these are the only obstructions for $F$ to be arc-extendible.

Suppose $F$ is a compact incompressible surface on the boundary of a compact, orientable, irreducible 3-manifold $M$. Let $F'$ be a component of $\partial M - \text{Int} F$. We say that $F$ is *arc-extendible* (in $F'$) if there is a properly embedded arc $\gamma$ in $F'$ such that $F \cup N(\gamma)$ is incompressible. In this case $\gamma$ is called an *extension arc* of $F$. We study the problem of which incompressible surfaces on the boundary $M$ are arc-extendible. This is useful in, for example, finding a sequence of mutually nonparallel incompressible surfaces in a 3-manifold.

Denote by $I$ the unit interval $[0, 1]$. We say that $M$ is a product $F \times I$ if there is a homeomorphism $\varphi : M \cong F \times I$ with $\varphi(F) = F \times 1$. Note that in this case $F' = \partial M - \text{Int} F$, and $F$ is not arc-extendible. A surface $F$ is *diskless* if it has no disk component. An incompressible surface with a disk component is always arc-extendible, unless the disk lies on a sphere component of $\partial M$. Thus to avoid trivial cases, we will only consider arc-extension of diskless surfaces.

**Theorem 1.** *Let $F$ be a diskless, compact, incompressible surface on the boundary of a compact, orientable, irreducible 3-manifold $M$, and let $F'$ be a non-annular component of $\partial M - \text{Int} F$. Then either $F$ is arc-extendible in $F'$, or $M$ is a product $F \times I$.*

The proof of the theorem involve some deep results about incompressible surfaces related to Dehn surgery and 2-handle additions. It breaks down into three cases. The

1991 *Mathematics Subject Classification.* Primary 57N10..

Freedman was partially supported by an NSF grant. Wu's research at MSRI was supported in part by NSF grant #DMS 9022140.

Typeset by $\mathcal{A}\mathcal{M}\mathcal{S}$-TEX

case that $F'$ is a thrice punctured sphere is treated in Theorem 4, which shows that if the surface obtained by gluing $F$ and $F'$ along one of the boundary curve of $F'$ is compressible for all the three boundary curves of $F'$, then $M$ must be a product. The second case is that $F'$ is parallel into $F$ (see below for definition). A similar result as above holds in this case. Theorem 9 shows that in the remaining case there is an arc $\gamma$ intersecting some circle $C$ in $F'$ at one point, so that all but at most three Dehn twists of $\gamma$ along $C$ are extension arcs of $F$. Moreover, in this case the extension arc $\gamma$ of $F$ can be chosen to have endpoints on any prescribed components of $\partial F'$. See Theorem 10 below.

Note that the irreducibility of $M$ is irrelevant to the compressibility of surfaces on $\partial M$. However, this does make the conclusion of the theorem simpler. If we drop this assumption from the theorem, the conclusion should be changed to "Either $F$ is arc-extendible in $F'$, or there is a component $F_0$ of $F$, and a homeomorphism $\varphi : M \cong F_0 \times I \# M'$ for some $M'$, such that $\varphi(F_0) = F_0 \times 1$, and $\varphi(F') = F_0 \times 0 \cup \partial F_0 \times I$."

Given a simple closed curve $\alpha$ on a surface $S$ on the boundary of $M$, we use $M[\alpha]$ to denote the manifold obtained by adding a 2-handle to $M$ along the curve $\alpha$. More explicitly, $M[\alpha]$ is the union of $M$ and a $D^2 \times I$, with the annulus $(\partial D^2) \times I$ glued to a regular neighborhood $N(\alpha)$ of $\alpha$ on $\partial M$. Use $S[\alpha]$ to denote the surface in $M[\alpha]$ corresponding to $S$, i.e. $S[\alpha] = (S - N(\alpha)) \cup (D^2 \times \partial I)$. The following two lemmas are very useful in dealing with incompressible surfaces. Various versions of Lemma 2 have been proved by Przytycki [Pr], Johannson [Jo], Jaco [Ja], and Scharlemann [Sch]. The lemma as stated is due to Casson and Gordon [CG].

**Lemma 2.** (The Handle Addition Lemma [CG].) *Let $\alpha$ be a simple closed curve on a surface $S$ on the boundary of an orientable irreducible 3-manifold $M$, such that $S$ is compressible and $S - \alpha$ is incompressible. Then $S[\alpha]$ is incompressible in $M[\alpha]$, and $M[\alpha]$ is irreducible.*

**Lemma 3.** (The Generalized Handle Addition Lemma.) *Let $S$ be a surface on the boundary of an orientable 3-manifold $M$, let $\gamma$ be a 1-manifold on $S$, and let $\alpha$ be a circle on $S$ disjoint from $\gamma$. Suppose $S - \gamma$ is compressible and $S - (\gamma \cup \alpha)$ is incompressible. If $D$ is a compressing disk of $S[\alpha]$ in $M[\alpha]$, then there is a compressing disk $D'$ of $S - \alpha$ in $M$ such that $\partial D' \cap \gamma \subset \partial D \cap \gamma$.*

*Proof.* This is essentially [Wu2, Theorem 1]. The theorem there stated that $\partial D' \cap \gamma$ has no more points than $\partial D \cap \gamma$, but the proof there gives the stronger conclusion that $\partial D' \cap \gamma \subset \partial D \cap \gamma$. □

We first study the case that the surface $F'$ in Theorem 1 is a thrice punctured sphere. Let $\alpha_1, \alpha_2, \alpha_3$ be the boundary curves of $F'$. Since $F'$ is a component of $\partial M - \text{Int} F$, we have $\alpha_i \subset \partial F$ for $i = 1, 2, 3$. Note that if $\text{Int} F \cup \text{Int} F' \cup \alpha_i$ is incompressible for some $i$, then for any essential arc $\gamma$ on $F'$ with $\partial \gamma \subset \alpha_i$, the surface $F \cup N(\gamma)$ is incompressible. Hence the following theorem proves Theorem 1 in the case that $F'$ is a twice punctured disk. However, it should be noticed that a similar statement is false if we drop the assumption that $F'$ is a sphere with three holes.

**Theorem 4.** *Let $F$ be a diskless compact incompressible surface on the boundary of a compact, orientable, irreducible 3-manifold $M$, and let $F'$ be a component of $\partial M - \text{Int} F$ which is a punctured sphere with $\partial F' = \alpha_1 \cup \alpha_2 \cup \alpha_3$. If $\text{Int} F \cup \text{Int} F' \cup \alpha_i$ is compressible for $i = 1, 2, 3$, then $M$ is a product $F \times I$.*

*Proof.* We fix some notation. Write $\widehat{F} = F \cup F'$. Denote by $\widehat{F}_i$ the surface obtained by gluing $\text{Int} F$ and $\text{Int} F'$ along $\alpha_i$, i.e. $\widehat{F}_i = \text{Int} F \cup \text{Int} F' \cup \alpha_i$. Similarly, write $\widehat{F}_{ij} = \text{Int} F \cup \text{Int} F' \cup \alpha_i \cup \alpha_j$.

First notice that $F'$ is incompressible. This is because each simple closed curve on $F'$ is isotopic to one of the $\alpha_i \subset F$, and because $F$ is incompressible and diskless. Since $\text{Int} F \cap \text{Int} F' = \emptyset$, the surface $\text{Int} F \cup \text{Int} F'$ is incompressible.

Let $M'$ be a maximal compression body of $\partial M$ in $M$. Then a surface on the boundary of $M$ is compressible in $M$ if and only if it is compressible in $M'$. Notice that if $M \neq M'$, then $M'$ is never a product $F \times I$, so if the theorem is true for $M'$, it is true for $M$. Hence after replacing $M$ by $M'$ if necessary, we may assume without loss of generality that $M$ is a compression body.

We claim that the curves $\alpha_1, \alpha_2, \alpha_3$ are mutually nonparallel on $\widehat{F}$, that is, no component of $F$ is an annulus with both boundary components on $F'$. If two curves $\alpha_1, \alpha_2$, say, are parallel on $\widehat{F}$, then the surface $\text{Int} F \cup \text{Int} F' = \widehat{F} - \alpha_1 \cup \alpha_2 \cup \alpha_3$ is incompressible if and only if $\widehat{F}_1 = \widehat{F} - \alpha_2 \cup \alpha_3$ is incompressible. However, by assumption $\widehat{F}_1$ is compressible, and we have shown that $\text{Int} F \cup \text{Int} F'$ is incompressible. Hence the claim follows.

Since $\widehat{F}_i$ is compressible, and $\widehat{F}_i - \alpha_i = \text{Int} F \cup \text{Int} F'$ is incompressible, we can apply the Handle Addition Lemma (Lemma 2) to $\widehat{F}_i$ and $\alpha_i$ to conclude that after adding a 2-handle along $\alpha_i$, the surface $\widehat{F}_i[\alpha_i]$ is incompressible in $M[\alpha_i]$, and $M[\alpha_i]$ is irreducible.

Consider the surface $\widehat{F}[\alpha_1]$. Notice that after adding the 2-handle, the surface $F'$ becomes an annulus on $\widehat{F}[\alpha_1]$ with boundary $\alpha_2 \cup \alpha_3$, so the two curves $\alpha_2, \alpha_3$ are parallel on $\widehat{F}[\alpha_1]$. Thus, $\widehat{F}_1[\alpha_1] = \widehat{F}[\alpha_1] - \alpha_2 \cup \alpha_3$ being incompressible in $M[\alpha_1]$ implies that $\widehat{F}[\alpha_1] - \alpha_2$ is incompressible in $M[\alpha_1]$. With the above notation, this says that $\widehat{F}_{13}[\alpha_1]$ is incompressible in $M[\alpha_1]$.

By assumption $\widehat{F}_3$ is compressible in $M$. Let $D$ be a compressing disk of $\widehat{F}_3$ in $M$. Then $\partial D$ is disjoint from $\alpha_1 \cup \alpha_2$, because $\partial D \subset \widehat{F}_3$. Also, $\partial D$ is not isotopic to $\alpha_1$ in $\widehat{F}_{13}$, otherwise $\alpha_1$ would bound a disk in $M$, contradicting the assumption that $F$ is diskless and incompressible. We have shown that $\widehat{F}_{13}[\alpha_1]$ is incompressible in $M[\alpha_1]$, so $D$ is not a compressing disk of $\widehat{F}_{13}[\alpha_1]$ in $M[\alpha_1]$, and hence $\partial D$ must bound a disk in $\widehat{F}_{13}[\alpha_1]$. This is true if and only if $\partial D$ is coplanar to $\alpha_1$ on $\widehat{F}_{13}$, that is, either $\partial D$ is parallel to $\alpha_1$, or it bounds a once punctured torus $T$ in $\widehat{F}_{13}$ which contains $\alpha_1$ as a nonseparating curve. The first possibility has been ruled out, so the second must be true. Let $\widehat{T}$ be the torus $T \cup D$. Since we have assumed above that $M$ is a compression body, either (i) $\widehat{T}$ is parallel to a boundary component of $M$, or (ii) $\widehat{T}$ bounds a solid torus.

If $\widehat{T}$ is parallel to a boundary component $T_0$ of $M$, then after adding the 2-handle, the surface $\widehat{T}[\alpha_1]$ becomes a sphere which separates $T_0$ from $\widehat{F}[\alpha_1]$, hence is a reducing sphere of $M[\alpha_1]$, which contradicts the irreducibility of $M[\alpha_1]$. Similarly, if $\widehat{T}$ bounds a solid torus $V$ but $\alpha_1$ is not a longitude of $V$, then after adding the 2-handle the manifold would have a lens space or $S^2 \times S^1$ summand, which again contradicts the irreducibility of $M[\alpha_1]$. (Note that $M[\alpha_1]$ cannot be a lens space because it has some boundary components.)

We have now shown that there is a compressing disk $D$ of $\widehat{F}_3$ in $M$ which cuts the manifold into two pieces, one of which is a solid torus $V$ which contains $\alpha_1$ as a longitude, but is disjoint from $\alpha_2$. Let $D_1$ be a meridian disk of $V$. Then $\partial D_1 \cap \alpha_1$ is a single point, and $\partial D_1$ is disjoint from $\alpha_2$ because $\partial V$ is disjoint from $\alpha_2$. Notice that $\partial D_1$ is not coplanar to $\alpha_2$, for if $\partial D_1$ were parallel to $\alpha_2$ then $\alpha_2$ would also intersect $\alpha_1$, and if $\partial D_1$ would bound a once punctured torus containing $\alpha_2$ then $\partial D_1$ would be a separating curve on $\partial M$, so it would intersect $\alpha_1$ in an even number of points, either case leading to a contradiction. Thus, after adding a 2-handle to $M$ along $\alpha_2$, the disk $D_1$ remains a compressing disk of $\widehat{F}[\alpha_2]$. Since the two curves $\alpha_1$ and $\alpha_3$ are parallel in $\widehat{F}[\alpha_2]$, and since $D_1$ intersects $\alpha_1$ in a single point, we can isotope $D_1$ to another disk $D_2$ in $M[\alpha_2]$ so that it intersects each of $\alpha_1$ and $\alpha_3$ in a single point. We are looking for such a disk in $M$; however $D_2$ is not necessary the one because it may intersect the attached 2-handle.

Recall that the surface $\widehat{F}_2$ is compressible, but the surface $\widehat{F}_2 - \alpha_2 = \text{Int} F \cup \text{Int} F'$ is incompressible. Hence we can apply the Generalized Handle Addition Lemma (Lemma 3, with $S = \widehat{F}$, $\gamma = \alpha_1 \cup \alpha_3$, and $\alpha = \alpha_2$) to conclude that there is also a compressing disk $D_3$ of $\widehat{F}$ in $M$, such that $\partial D_3$ is disjoint from $\alpha_2$, and $\partial D_3 \cap (\alpha_1 \cup \alpha_3)$ is a subset of $\partial D_2 \cap (\alpha_1 \cup \alpha_3)$.

The set $\partial D_3 \cap (\alpha_1 \cup \alpha_3)$ is nonempty, otherwise, since $\partial D_3$ is also disjoint from $\alpha_2$, $D_3$ would be a compressing disk of $\text{Int} F \cup \text{Int} F'$, contradicting the incompressibility of $\text{Int} F \cup \text{Int} F'$. Since $\alpha_1 \cup \alpha_2 \cup \alpha_3$ is separating on $\widehat{F}$, the curve $\partial D_3$ can not intersect $\alpha_1 \cup \alpha_2 \cup \alpha_3$ at a single point. It follows that $\partial D_3 \cap (\alpha_1 \cup \alpha_3) = \partial D_2 \cap (\alpha_1 \cup \alpha_3)$, that is, $\partial D_3$ intersects each of $\alpha_1, \alpha_3$ in a single point. Such a disk is called a *bigon*.

Denote by $D_{13}$ the bigon $D_3$ above. Interchanging the rules of $\alpha_1$ and $\alpha_2$ in the above argument, we get another compressing disk $D_{23}$ of $\widehat{F}$ in $M$, which is disjoint from $\alpha_1$, and intersects each of $\alpha_2, \alpha_3$ in a single point. By a simple innermost circle outermost arc argument, we can isotope $D_{13}$ so that it is disjoint from $D_{23}$, and still has the same number of intersection points with each $\alpha_i$. Cutting $M$ along $D_{13} \cup D_{23}$, we get a submanifold $M'$ of $M$, in which the surface $F'$ becomes a disk $\widetilde{F}' \subset F'$, and the surface $F$ becomes a surface $\widetilde{F} \subset F$. It is clear that one boundary component $C$ of $\widetilde{F}$ bounds a disk on $\partial M'$, namely the union of $\widetilde{F}'$ and the two copies of $D_{13} \cup D_{23}$. Since $F$ is incompressible, this curve $C$ bounds a disk in $F$, so $\widetilde{F}$ must be a disk. These disks together form a sphere boundary component of $M'$. Since $M$ is irreducible, $M'$ must be a 3-ball, so it is a product $\widetilde{F} \times I$. Gluing back along $D_{13}$ and $D_{23}$, we see that $M$ is a product $F \times I$. This completes the proof of Theorem 4. $\square$

Below, $F, F'$ and $M$ will be as in Theorem 1. Using Theorem 4 we may assume that $F'$ is not a thrice punctured sphere. A curve $C'$ on $F'$ is *$\partial$-nonseparating* if (i) $C'$ is not parallel to a boundary curve on $F'$, and (ii) there is a proper arc $\gamma$ in $F'$ intersecting $C'$ in a single point. A sub-surface $G'$ of $F'$ is *parallel into $F$* if there is a product $G' \times I \subset M$ such that $G' = G' \times 0$, and $G' \times 1 \subset F$. Similarly, a curve $C'$ on $F'$ is *parallel into $F$* if there is an embedded annulus $A \subset M$ with $\partial A = C' \cup C$, where $C \subset F$.

**Lemma 5.** *If $F'$ is compressible, then there is a $\partial$-nonseparating curve $C'$ on $F'$ which is not parallel into $F$.*

*Proof.* Let $D$ be a compressing disk of $F'$. If $\partial D$ is non-separating on $F'$, let $C'$ be a curve in $F'$ that intersects $\partial D$ in one point. Then $C'$ is nonseparating, hence $\partial$-nonseparating on $F'$. We want to show that $C'$ is not parallel into $F$. Otherwise, let $A$ be an annulus with $\partial A = C' \cup C$, where $C \subset F$. Then $A \cap D$ is a proper 1-manifold on $D$. But $\partial(A \cap D) = (\partial A) \cap \partial D$ is a single point, which is absurd. Hence $C'$ is the curve required.

Now assume that $\partial D$ is separating on $F'$, cutting $F'$ into $F'_1$ and $F'_2$. Choose a simple loop $C_i$ on $F'_i$ as follows. If $F'_i$ is nonplanar, then there are a pair of non-separating curves intersecting each other in one point, at least one of which is not null-homologous in $M$. Choose this one as $C_i$. If $F'_i$ is planar, choose $C_i$ to be isotopic to a boundary curve of $F'$. Note that since $F$ is incompressible and diskless, $C_i$ is not null-homotopic in $M$. Also notice that in both cases there is a properly embedded arc $\gamma$ on one of the $F'_i$ which intersects $C_1 \cup C_2$ in one point.

Now choose a band $B = I \times I$ on $F'$ such that $B \cap \partial D = I \times \frac{1}{2}$, $B \cap C_1 = I \times 0$, $B \cap C_2 = I \times 1$, and $B$ is disjoint from the arc $\gamma$ above. Such band exists because $\gamma$ is a nonseparating arc on $F'_i$. Let $C'$ be the band sum of $C_1$ and $C_2$, that is, $C' = (C_1 \cup C_2 - I \times \{0,1\}) \cup (\{0,1\} \times I)$. Then $C'$ intersects $\gamma$ in one point. Since $C'$ intersects $\partial D$ essentially in two points, it is not parallel to any boundary component on $F'$. Therefore $C'$ is $\partial$-nonseparating.

We want to show that $C'$ is not parallel into $F$. Using the property that $C_i$ are not null-homotopic in $M$, one can show by an innermost circle argument that $C'$ is not null-homotopic in $M$. Now suppose that there is an annulus $A$ in $M$ with $\partial A = C' \cup C$, where $C \subset F$. Since $C'$ is not null-homotopic in $M$, $A$ is incompressible in $M$. By surgery along an innermost circle of $D \cap A$ one can eliminate all circle intersections of $A \cap D$. Since $\partial(A \cap D)$ consists of two points, $A \cap D$ is a single arc, which has endpoints on the same component of $\partial A$, hence it cuts off a disk $D'$ from $A$. Assume without loss of generality that $D' \cap F'$ is on $F'_1$. Let $D''$ be the disk on $D$ bounded by $(A \cap D) \cup (B \cap D)$, and let $B_1 = B \cap F'_1$. Then $D' \cup D'' \cup B_1$ is a disk with boundary $C_1$, which contradicts the fact that $C_1$ is not null-homotopic in $M$. Therefore, $C'$ is not parallel into $F$. □

**Lemma 6.** *Suppose $F'$ is incompressible, and is not a thrice punctured sphere. Then either (i) there is a $\partial$-nonseparating curve $C'$ on $F'$ which is not parallel into $F$, or (ii) $F'$ is parallel into $F$.*

*Proof.* Since $F'$ is not a thrice punctured sphere, one can easily find a $\partial$-nonseparating curve $\alpha_0$ on $F'$. Assume that (i) is not true, so all $\partial$-nonseparating curves are parallel into $F$. We want to show that $F'$ is parallel into $F$.

Since $\alpha_0$ is parallel into $F$, the annulus $N(\alpha_0)$ is also parallel into $F$. It is an incompressible annulus because $\alpha_0$ is essential on $F'$ and $F'$ is incompressible. Among all connected incompressible surfaces in Int$F'$ which contain $\alpha_0$ and are parallel into $F$, choose $G'$ such that the complexity $(\chi(G'), |\partial G'|)$ is minimal in the lexical-graphic order, where $\chi(G')$ is the Euler characteristic of $G'$, and $|\partial G'|$ is the number of boundary components of $G'$. All incompressible sub-surfaces of $F'$ have Euler characteristics bounded below by $\chi(F')$, hence such $G'$ does exist.

If all boundary components of $G'$ are parallel to some boundary components on $F'$, then either $G'$ is contained in a collar of $\partial F'$, or $F' - \text{Int}G' = \partial F' \times I$. The first case does not happen because $G'$ contains the $\partial$-nonseparating curve $\alpha_0$, which by definition is not parallel to any boundary curve on $F'$. In the second case $F'$ is isotopic to $G'$, so it is parallel into $F$, and we are done. Hence we may assume that some boundary curve $\beta$ of $G'$ is not parallel to any boundary curve on $F'$.

We want to find a $\partial$-nonseparating curve $\alpha'$ which intersects $\beta$ essentially in one or two points. If $\beta$ is nonseparating on $F'$, choose $\alpha'$ to be any curve on $F'$ that intersects $\beta$ in a single point. Then $\alpha'$ is nonseparating, hence $\partial$-nonseparating on $F'$.

If $\beta$ separates $F'$ into $F'_1$ and $F'_2$, choose an essential arc $\alpha'_i$ on $F'_i$ with $\partial\alpha'_1 = \partial\alpha'_2 \subset \beta$. Moreover, if $F'_i$ is nonplanar, choose $\alpha'_i$ to be nonseparating on $F'_i$. Then $\alpha' = \alpha'_1 \cup \alpha'_2$ is $\partial$-nonseparating, and intersects $\beta$ essentially in two points, as required.

Isotope $\alpha'$ so that it intersects $\partial G'$ minimally. The geometric intersection number between $\alpha'$ and $\beta$ is 1 or 2, so $\alpha' \cap \partial G' \neq \emptyset$. Since $\alpha'$ is $\partial$-nonseparating, by our assumption above it is parallel into $F$, so there is an annulus $A$ with $\partial A = \alpha' \cup \alpha$, where $\alpha \subset F$. Isotope $A$ rel $\alpha'$ so that it intersects $(\partial G') \times I$ minimally. Since $G'$ is incompressible, $(\partial G') \times I$ consists of incompressible annuli in $M$, hence $A \cap ((\partial G') \times I)$ has no trivial circles. Since $F$ and $F'$ are also incompressible, one can show that $A \cap ((\partial G') \times I)$ has no trivial arcs on $A$ either. Therefore $A \cap ((\partial G') \times I)$ consists of vertical arcs $(\alpha' \cap \partial G') \times I$. These arcs cut $A$ into several squares $\alpha'_i \times I$, where each $\alpha'_i$ is the closure of a component of $\alpha' - \partial G'$. Choose $i$ so that $\alpha'_i$ lies outside of $G'$. Let $H$ be the component of $F' - \text{Int} G'$ that contains $\alpha'_i$. Then $G'' = G' \cup N(\alpha'_i)$ is a surface parallel into $F$, and $\chi(G'') = \chi(G') - 1$. The arc $\alpha'_i$ is essential on $H$, so the only case that some boundary component $\gamma$ of $G''$ bounds a disk on $F'$ is when $H$ is an annulus, and $\gamma$ is the boundary of the disk obtained by cutting $H$ along $\alpha'_i$. Since $F$ and $F'$ are incompressible and $M$ is irreducible, both ends of the annulus $\gamma \times I \subset G'' \times I \subset M$ bound disks on $F \cup F'$, which together with $\gamma \times I$ bounds a 3-ball in $M$. It follows that $G' \cup H$ is parallel into $F$. Since $G' \cup H$ has the same Euler characteristic as $G'$ but fewer number of boundary components, this contradicts the choice of $G'$. Therefore $\partial G''$ consists of essential curves on $F'$. Since $F'$ is incompressible, $G''$ is also incompressible. Since $\chi(G'') < \chi(G')$, this again contradicts the choice of $G'$. $\square$

Given a simple closed curve $\alpha$ and a proper arc $\gamma$ on $F'$, denote by $\tau^n_\alpha \gamma$ the curve obtained from $\gamma$ by Dehn twist $n$ times along $\alpha$, and by $N(\tau^n_\alpha \gamma)$ a regular neighborhood of $\tau^n_\alpha \gamma$ on $\partial M$. Suppose $T$ is a fixed torus boundary component of a 3-manifold $M$. Denote by $M(r)$ the manifold obtained by Dehn filling on $T$ along a slope $r$ on $T$, that is $M(r)$ is obtained by gluing a solid torus $V$ to $M$ along $T$ so that the curve $r$ on $T$ bounds a meridian disk on $V$. Denote by $\Delta(r_1, r_2)$ the minimal geometric intersection number between two slopes $r_1, r_2$. The following two theorems will be used in the proof of Theorem 9, which proves Theorem 1 in the case that $F'$ contains a $\partial$-nonseparating curve which is not parallel into $F$.

**Lemma 7.** ([Wu2], Theorem 1) *Let $T$ be a torus component on the boundary of a 3-manifold $M$, and let $S$ be an incompressible surface on $\partial M - T$. Suppose there is no incompressible annulus in $M$ with one boundary component on each of $S$ and $T$. If $S$ is compressible in $M(r_1)$ and $M(r_2)$, then $\Delta(r_1, r_2) \leq 1$. In particular, $S$ is incompressible in all but at most three $M(r)$.* $\square$

**Lemma 8.** ([CGLS], Theorem 2.4.3) *Let $T, S, M$ be as in Lemma 7, and assume further that $M$ is irreducible. Suppose that there is an incompressible annulus $A$ in $M$ with one boundary component on $S$ and the other a curve $r_0$ on $T$. Then either $S$ is a torus and $M = S \times I$, or $S$ remains incompressible in all $M(r)$ with $\Delta(r, r_0) > 1$.* □

**Theorem 9.** *Let $\alpha$ be a $\partial$-nonseparating curve on $F'$ which is not parallel into $F$, and let $\gamma$ be a proper arc on $F'$ intersecting $\alpha$ in one point. Then $F_n = F \cup N(\tau_\alpha^n \gamma)$ is incompressible for all but at most three consecutive $n$'s.*

*Proof.* Let $K$ be the knot obtained by pushing $\alpha$ slightly into $M$. There is an embedded annulus $A_0$ in $M$ with $\partial A_0 = \alpha \cup K$. Consider the manifold $M_K = M - \text{Int} N(K)$, where $N(K)$ is a regular neighborhood of $K$ in $M$. Let $T$ be the torus $\partial N(K)$, and let $(m, l)$ be the meridian-longitude pair on $T$ such that $l = A_0 \cap T$. Denote by $M_K(p/q)$ the manifold obtained by Dehn filling on $T$ along the slope $pm + ql$. The Dehn twist $\tau_\alpha^{-n}$ on $F'$ extends to a Dehn twist of $M_K$ along the annulus $A = A_0 \cap M_K$, which sends the meridian slope $m$ of $T$ to the slope $m - nl$, so it extends to a homeomorphism $\varphi_n : M = M_K(1/0) \cong M_K(-1/n)$, which maps the curve $\tau_\alpha^n \gamma$ to the curve $\gamma$, and hence the surface $F_n$ to the surface $F_0 = F \cup N(\gamma)$. It follows that $\varphi_n$ is a homeomorphism of pairs

$$\varphi_n : (M, F_n) \to (M_K(-1/n), F_0).$$

Therefore to prove the theorem we need only show that for all but at most three consecutive integers $n$, the surface $F_0$ is incompressible in $M_K(-1/n)$.

CLAIM 1. *$T = \partial N(K)$ is incompressible in $M_K$, and $M_K$ is irreducible.*

If $D$ is a compressing disk of $T$ in $M_K$, then $\partial D$ must intersect the meridian $m$ of $K$ in one point, because otherwise after the trivial Dehn filling, $M = M_K(1/0)$ would contain a lens space or $S^2 \times S^1$ summand, contradicting the irreducibility of $M$. It follows that $K$, and hence $\alpha$, bounds a disk in $M$. In this case $\alpha$ is parallel to a trivial curve on $F$, which contradicts the assumption that $\alpha$ is not parallel into $F$. Similarly, if $M_K$ is reducible, then since $M$ is irreducible, $K$ is contained in a ball in $M$, so $\alpha$ would be null-homotopic. Using Dehn's Lemma, we see that $\alpha$ bounds a disk in $M$, hence is parallel to a trivial circle in $F$, contradicting the assumption that $\alpha$ is not parallel into $F$.

CLAIM 2. *$F_0$ is incompressible in $M_K$.*

Recall that $A$ denotes the annulus $A_0 \cap M_K$. Since $\alpha$ intersects $\gamma$ in a single point, $A \cap F_0$ is a single arc $C$ on the boundary curve $\alpha$ of $A$. Let $D$ be a compressing disk of $F_0$, chosen so that $|D \cap A|$, the number of components in $D \cap A$, is minimal. After disk swapping along disks on $A$ bounded by innermost circles, we can assume that no component of $D \cap A$ is a trivial circle on $A$. Since $T$ is incompressible by Claim 1, the annulus $A$ is also incompressible, so $D \cap A$ contains no essential circle component on $A$ either. Hence $D \cap A$ consists of arcs only. If some arc $e$ of $D \cap A$ is parallel to a sub-arc on $C = A \cap F_0$, then after boundary compressing $D$ along a disk $\Delta$ cut off by an outermost such arc we will get two disks $D_1, D_2$ with boundary on $F_0$, at least one of which has boundary an essential curve on $F_0$, hence is a compressing disk of $F_0$. Since $|D_i \cap A| < |D \cap A|$, this contradicts the minimality of $|D \cap A|$. Therefore, all arcs of $D \cap A$ are essential relative to $C$, in the sense that it is not parallel to an arc on $C$. See Figure 1(a). Notice that $|D \cap A| \neq 0$, otherwise $D$ would be a compressing disk of $F$, contradicting the incompressibility of $F$.

Consider an outermost disk $\Delta$ on $D$, as shown in Figure 1(b). Then $\partial \Delta$ consists of two arcs $e_1, e_2$, where $e_1$ is an arc on $A$ which is essential relative to $C$, and $e_2$ is an arc on $F_0$ with interior disjoint from $C$. Thus $e_2 \cap N(\gamma)$ consists of two arcs $e_2', e_2''$. Let $t_1$ be the subarc of $C$ connecting the two ends of $e_2' \cup e_2''$ on $C$, and let $t_2$ be the subarc on $\partial N(\gamma)$ connecting the other two ends of $e_2' \cup e_2''$. Then $e_2' \cup t_1 \cup e_2'' \cup t_2$ bounds a disk $\Delta'$ on $N(\gamma)$. Now $A' = \Delta \cup \Delta'$ is an annulus in $M$, with one boundary component $e_1 \cup t_1$ an essential circle on $A$, which is parallel to $\alpha$, and the other component $e_2 \cup t_2$ a curve on $F$. This contradicts the assumption that $\alpha$ is not parallel into $F$.

Figure 1

CLAIM 3. *There is no incompressible annulus $P$ in $M_K$ with one boundary component $C_1$ on $F_0$ and the other component $C_2$ a curve on $T$ which is disjoint from $l = A \cap T$.*

The proof is similar to that of Claim 2. Choose $P$ so that $|P \cap A|$ is minimal. Using the fact that $P$ is incompressible, one can show as above that $P \cap A$ has no trivial circle component. Note that since $C_2$ is disjoint from $l$, $P \cap A$ has no arc component with endpoints on $l = A \cap T$. If $P \cap A$ had some essential circle component, choose such a component $t$ which is closest to $l$ on $A$. By cutting and pasting along the annulus on $A$ bounded by $t \cup l$, one would get another incompressible annulus $P'$ which has fewer intersection components with $A$. As in the proof of Claim 2 one can eliminate all arc components of $P \cap A$ which are inessential relative to $C = A \cap F_0$. Hence $P \cap A$ consists of arcs with ends on $C$ and are essential relative to $C$, as shown in Figure 1(a). Also, since $P$ is disjoint from $l$, $P \cap A$ are inessential arcs on $P$. Now one can use a disk $\Delta$ cut off by an outermost arc on $P$, proceed as in the proof of Claim 2 to get an annulus with one boundary on $\alpha$ and the other on $F$, and get a contradiction. Finally, if $P \cap A = \emptyset$ then $P$ extends to an annulus with one boundary on $\alpha$ and the other on $F$, contradicting the assumption that $\alpha$ is not parallel into $F$. This completes the proof of Claim 3.

We now continue with the proof of Theorem 9. We have shown that $F_0$ is incompressible in $M_K$. If there is no essential annulus in $M_K$ with one boundary component on each of $F_0$ and $T$, then by Lemma 7 we know that $F_0$ is incompressible in $M_K(r)$ for all but at most three slopes $r$ with mutual intersection number 1. In particular, it is incompressible in $M_K(-1/n)$ for all but at most two consecutive $n$'s, so the theorem follows. Now suppose there is an essential annulus $P$ in $M_K$ with one boundary component on $F_0$ and the other a curve $r_0$ on $T$. Since $F_0$ is not a closed surface, it is not a torus. Hence by Lemma 8, $F_0$ remains incompressible in $M_K(-1/n)$ unless $\Delta(-1/n, r_0) \leq 1$. By Claim 3, $r_0$ is not the longitude slope $0/1$, therefore, $\Delta(-1/n, r_0) \leq 1$ holds for at most three consecutive integers $n$. This completes the proof of Theorem 9. □

*Proof of Theorem 1.* By Theorem 4, Lemmas 5 and 6, and Theorem 9, we can now assume that $F'$ is incompressible and is parallel into $F$. We want to show that either $F$ is arc-extendible in $F'$, or $M$ is a product $F \times I$. As in the proof of Theorem 4, we may assume without loss of generality that $M$ is a compression body, so all closed incompressible surfaces of $M$ are boundary parallel. Let $\alpha_1, \ldots, \alpha_k$ be the boundary curves of $F'$. Let $F' \times I$ be a product in $M$ such that $F' = F' \times 0$ and $F' \times 1 \subset F$. Write $\alpha_i^1 = \alpha_i \times 1$, which is a curve on $F$ isotopic to $\alpha_i$ in $M$.

We have assumed above that $F'$ is incompressible in $M$, so $\operatorname{Int} F \cup \operatorname{Int} F'$ is incompressible in $M$. Write $\widehat{F}_i = \operatorname{Int} F \cup \operatorname{Int} F' \cup \alpha_i$. If $\widehat{F}_i$ is incompressible for some $i$, then $F \cup N(\gamma)$ is incompressible for any essential arc $\gamma$ in $F'$ with endpoints on $\alpha_i$, and we are done. (Such an arc exists because $F'$ is not an annulus or disk.) So assume that $\widehat{F}_i$ is compressible for all $i$. By the Handle Addition Lemma (Lemma 2), after adding a 2-handle to $M$ along $\alpha_i$, the surface $\widehat{F}_i[\alpha_i]$ is incompressible, and $M[\alpha_i]$ is irreducible. Notice that since $F'$ is incompressible, the curve $\alpha_i^1 = \alpha_i \times 1$ in $F$ is essential on $F$. But after adding the 2-handle, $\alpha_i^1$ bounds a disk in $M[\alpha_i]$, so it must also bound a disk on $\widehat{F}_i[\alpha_i]$ because $\widehat{F}_i[\alpha_i]$ is incompressible. By definition $\widehat{F}_i[\alpha_i]$ is obtained from $(\mathrm{Int} F \cup \mathrm{Int} F') - \mathrm{Int} N(\alpha_i)$ by capping off the two copies of $\alpha_i$ with disks, hence $\alpha_i^1 \cup \alpha_i$ bounds an annulus $A_i$ on $F_i$. Denote by $A'_i$ the annulus $\alpha_i \times I \subset F' \times I \subset M$. Then $T_i = A_i \cup A'_i$ is a torus in $M$. Since we have assumed above that $M$ is a compression body, either $T_i$ bounds a solid torus $V_i$, or it is parallel to some torus component of $\partial M$. However, since $M[\alpha_i]$ is irreducible, one can show as in the proof of Theorem 4 that $V_i$ is a solid torus, and $\alpha_i$ is a longitude of $V_i$. This is true for all $i$. It is now easy to see that $M$ is a product $F \times I$. □

The following theorem supplements Theorem 1. It says that in most case there are extension arcs with endpoints on any prescribed boundary compponents of $F'$.

**Theorem 10.** *Let $F, F', M$ be as in Theorem 1. Suppose $M$ is not a product $F \times I$, and suppose $F'$ is not parallel into $F$ and is not a thrice punctured sphere. Then it contains an extension arc $\gamma$ of $F$ with endpoints on any prescribed components of $\partial F'$.*

*Proof.* If $F'$ is nonplanar, then by the proof of Lemmas 5 and 6, there is a $\partial$-nonseparating circle $\alpha$ (denoted by $C'$ there) on $F'$ which is not parallel into $F$, and is actually nonseparating on $F'$. Hence given any boundary components $\partial_1, \partial_2$ of $F'$, (possibly $\partial_1 = \partial_2$), there is an arc $\gamma$ with endpoints on $\partial_1$ and $\partial_2$, intersecting $\alpha$ in one point. By Theorem 9, for all but at most three integers $n$, the arc $\gamma_n = \tau_\alpha^n \gamma$ is an extension arc of $F$.

Now suppose $F'$ is planar with $|\partial F'| \geq 4$. First assume that $\partial_1, \partial_2$ are distinct boundary components of $F'$. By the proof of Lemmas 5 and 6, the curve $\alpha$ is a band sum of two boundary components of $F'$. From the proofs one can see that we can always choose $\alpha$ to be the band sum of $\partial_1$ and $\partial_3$, with $\partial_3 \neq \partial_1, \partial_2$. Hence there is an arc $\gamma$ from $\partial_1$ to $\partial_2$ intersecting $\alpha$ in one point. We can then apply Theorem 9 to get an extension arc $\gamma_n$ with one endpoint on each of $\partial_1$ and $\partial_2$.

We now proceed to find an extension arc in $F'$ with boundary on the same component $\partial_1$ of $\partial F'$. By the proof of Lemmas 5 and 6, we can choose the curve $\alpha$ above to be the band sum of of $\partial_2$ and $\partial_3$, with $\partial_1 \neq \partial_2, \partial_3$. Recall that $\alpha$ is not parallel into $F$. Choose an arc $\gamma$ as follows. Let $\partial'_2$ be a curve on $F'$ parallel to $\partial_2$, let $\gamma'$ be an arc connecting $\partial'_2$ to $\partial_1$ intersecting $\alpha$ in one point, and let $Q$ be the sub-surface $N(\gamma' \cup \partial'_2)$ of $F'$. Then $\gamma$ is the closure of the arc component of $\partial Q \cap \mathrm{Int} F'$, that is, $\gamma$ is the arc component of the frontier of $Q$ in $F'$, see Figure 2 below. Consider the surface $F_0 = F \cup N(\gamma)$, and observe that $F_0$ is isotopic to the surface $F \cup Q$. After Dehn twist along $\alpha$, it is isotopic to the surface $F \cup N(\tau_\alpha^n \gamma)$; hence to show that all but at most three $\tau_\alpha^n \gamma$ are extension arcs of $F$ in $F'$, we need only show that $F \cup Q$ is incompressible after all but at most three Dehn twist along $\alpha$. Since $F \cup Q$ intersects $\alpha$ in a single arc, the argument in the proof of Theorem 9 is still valid, with the following easy modifications. We use the notations in that proof.

Figure 2

The proof of Claim 2 needs the following modifications. (i) The arc $e_2$ on the boundary of the outermost disk $\Delta$ may be on $Q$. In this case, notice that the other arc $e_1$ on $\partial D$ is isotopic to an arc $\alpha_1$ on $\alpha$, and $e_2 \cup \alpha_1$ is isotopic in $F'$ to the curve $\partial_3$, so the fact that $e_1 \cup e_2$ bounds a disk $\Delta$ would imply that $\partial_3$ bounds a disk. Since $\partial_3$ is also on $\partial F$, this contradicts the fact that $F$ is incompressible and diskless. (ii) The compressing disk $D$ of $F \cup Q$ could be disjoint from the annulus $A$. But since $F$ is incompressible, this would imply that $\partial D$ lies on $Q$, hence is isotopic to $\partial_2$, which would imply that $\partial_2$ bounds a disk, again contradicting the assumption that $F$ is incompressible and diskless.

The proof of Claim 3 applies to show that the annulus $P$ there can be modified to be disjoint from the annulus $A$. Then notice that the component of $\partial P$ on $F \cup Q$ is either in $F$, or in $Q$ and hence parallel to $\partial_2$. Since $\partial_2 \subset F$, in either case $P$ can be extended to an annulus with one boundary component on $\alpha$ and the other on $F$, which contradicts the assumption that $\alpha$ is not parallel into $F$.

The rest part of the proof of Theorem 9 follows verbatim to show that $F \cup Q$ is incompressible after all but at most three Dehn twist along $\alpha$.  □

*Remark.* Theorem 10 is not true if $F'$ is a thrice punctured sphere.

MICHAEL FREEDMAN DEPARTMENT OF MATHEMATICS, UC SAN DIEGO, LA JOLLA, CA 92093
*E-mail address*: `freedman@math.ucsd.edu`

HUGH HOWARDS DEPARTMENT OF MATHEMATICS, UC SAN DIEGO, LA JOLLA, CA 92093
*E-mail address*: `howards@math.ucsd.edu`

YING-QING WU DEPARTMENT OF MATHEMATICS, UNIVERSITY OF IOWA, IOWA CITY, IA 52242; and
MSRI, 1000 CENTENNIAL DRIVE, BERKELEY, CA 94720-5070
*E-mail address*: `wu@math.uiowa.edu`